\documentclass[12pt]{article}
\usepackage{amssymb,amsmath, mathrsfs,textcomp} 

\usepackage{graphicx}        

\newcommand{\pr}{{\bf P}}                     
\newcommand{\exn}{{\bf E\,}}                    
\newcommand{\ind}{{\bf 1}}                     
\newcommand{\ep}{\varepsilon}            
\newcommand{\la}{\lambda}            
\newcommand{\bla}{\boldsymbol \lambda}            
\newcommand{\La}{\Lambda}            
\newcommand{\al}{\alpha}                 
\newcommand{\de}{\delta}                 
\newcommand{\bal}{\boldsymbol \alpha}                 
\newcommand{\wbal}{\widehat \bal}                       
\newcommand{\bz}{\boldsymbol \zeta}                 
\newcommand{\by}{\boldsymbol y}                 
\newcommand{\bnu}{\boldsymbol{\nu}}                
\newcommand{\be}{\boldsymbol e}                 

\newcommand{\bX}{\boldsymbol X}                 

\newcommand{\bZ}{\boldsymbol Z}                 
\newcommand{\cA}{\mathcal{A}}             
\newcommand{\deq}{\stackrel{d}{=}}      

\newcommand{\bx}{\boldsymbol x}                 

\newcommand{\R}{\mathbb{R}}             
\newcommand{\Z}{\mathbb{Z}}             

\newtheorem{theo}{Theorem}

\newtheorem{rema}{Remark}
\newtheorem{lemo}{Lemma}

\newcommand{\zhe}[6]{#1 #2. {\em #3}, {\bf #5} (#4), #6.}
\newcommand{\kn}[4]{#1 {\em #2}. #3 (#4).}

\topmargin= - 8  true mm
\title{On large deviation probabilities for self-normalized sums of random variables\footnote{Research supported by the University of Melbourne Faculty of Science Research Grant Support Scheme.}}
\date{}
\author{Konstantin Borovkov\footnote{School of Mathematics and Statistics, The University of Melbourne, Parkville 3010, Australia; e-mail: borovkov@unimelb.edu.au.}} 

\begin{document}

\maketitle

\begin{abstract}
We reduced the large deviation problem for a self-normalized random walk to  one for an auxiliary usual bivariate random walk. This enabled us to prove the classical theorem for self-normalized walks by Q.-M.~Shao~(1997) under slightly more general conditions and, moreover, to provide a graphical interpretation for the emerging limit in terms of the rate function for the bivariate problem. Furthermore, using this approach, we obtained exact (rather than just logarithmic) large deviation asymptotics for the probabilities of interest. Extensions to more general self-normalizing setups including the multivariate case were discussed.

\medskip
{\em AMS Mathematics Subject Classification:}   60F10, 60G50.

\medskip
{\em Key words and phrases:} self-normalized partial sums, random walk, multivariate large deviations, exact large deviations asymptotics.

\end{abstract}

\section{Introduction}
\label{Sect_1}

Let $X, X_1, X_2,\ldots$ be a sequence of  independent and identically
distributed (i.i.d.) non-degenerate  random variables, and set 
\[
S_n:= \sum_{j=1}^n X_j, \quad n=1,2,\ldots
\]
The classical large deviations (LD) theory for the sequence of partial sums $\{S_n\}_{n\ge 1}$ establishes, under suitable moment or more specific distributional tail conditions on~$X$, the asymptotic behavior of the probabilities of the form $\pr (S_n >z_n)$ as $n\to\infty$ when the sequence $\{z_n\}$ is such that these probabilities are vanishing. This area of probability theory nicely complements the fundamental field of limit theorems on weak convergence of the distributions of the normalized sums $S_n$ and has numerous important uses in both theoretical and applied problems. Recent accounts of the LD theory in such and more general settings can be found in monographs~\cite{Bo20} and~\cite{BoBo08}, treating the case of light- and heavy-tailed distributions, respectively, and in the references therein. For the LD theory expositions in a broader context, we refer the reader to monographs~\cite{DeZe98} and~\cite{Pu01}. Note that most of the publications on the LD theory (including the two just mentioned monographs) deal with the ``crude" (or logarithmic) LD asymptotics; for ``exact LD asymptotics" (as presented in Section~\ref{Sect_3} of this paper), one can find an introductory level exposition in Section~9.3 of~\cite{Bo13} (in the univariate case only) and advanced treatments in~\cite{BoMo92}, \cite{Bo20}.

Motivated by the Griffin--Kuelbs~\cite{GrKu89} self-normalized law of the iterated logarithm for all distributions in the domains of attraction of stable laws, the seminal   paper~\cite{Sh97}  aimed to establish   a ``self-normalized LD theory'' for arbitrary random variables without any moment conditions (for more limit theorems in self-normalized settings, see more recent monograph~\cite{deLaSh08} and survey paper~\cite{ShWa13}). For $p\ge 0,$ set
\[
T_n=T_{n,p}:=\sum_{j=1}^n |X_j|^p, \quad n=1,2,\ldots
\]
The first key result from   paper~\cite{Sh97} is stated in its Theorems~1.2 (which includes Theorem~1.1 from the same paper as a special case when $p=2$) and claims the following about the LDs of the self-normalized sums
\begin{align}
\label{Wnp}
W_{n,p}:=\left\{
\begin{array}{ll}
\displaystyle \frac{S_n }{n^{1-1/p}T_n^{1/p}} & \mbox{if } T_n>0,
\\
\infty &  \mbox{if } T_n=0.
\end{array}
\right.
\end{align}
The convention about the infinite value of $W_{n,p}$ on the event  $
\{T_n=0\}=\{X_1=\cdots = X_n=0\}$
was stated in  Theorem~1.1 in~\cite{Sh97}; we  use   it in this paper as well.

\begin{theo}
\label{thm_Shao}
Let $p>1.$ Assume that either  $\exn X\ge 0$ or $\exn |X |^p =\infty.$ Then, for
\begin{align}
\label{z*}
z> z^*:=\left\{
\begin{array}{ll}
\exn X (\exn |X|^p)^{-1/p} & \mbox{if $\exn |X|^p<\infty$,}
\\
0 & \mbox{otherwise,}
\end{array}
\right.
\end{align}
there exists the limit
\begin{align}
\label{Shao_thm}
\lim_{n\to \infty } \frac1n & \ln\pr ( W_{n,p}
   \ge z  )=J_z,
\end{align}
where
\begin{align}
\label{Shao_thm_1}
 J_z:=
 \sup_{c\ge 0} \inf_{t\ge 0} \ln
 \exn \exp \bigl\{
 t (cX - zp^{-1}(|X|^p + (p-1) c^{p/p-1}))\bigr\}.
\end{align}
\end{theo}

\begin{rema}\label{rem_1}
\rm
By H\"older's inequality,   $\big|n^{-1}\sum_{j=1}^n a_j\big|\le n^{-1/p}\big(\sum_{j=1}^n |a_j|^p\big)^{1/p}$ for any real $a_j$'s and $p\ge 1,$ so that one always has
\begin{align}
\label{S<T}
{|S_n |}/(n^{1-1/p}T_n^{1/p})\le 1\quad \mbox{provided that $T_n>0.$}
\end{align}
Therefore, according to the above-mentioned convention in the second line of~\eqref{Wnp}, for $z>1,$ both sides in~\eqref{Shao_thm} are equal to $\ln \pr (X=0)\ge -\infty.$ For $z\in (z^*,1),$ one has $J_z\in (-\infty, 0),$ as will be seen from Lemma~\ref{prop} below.
\end{rema}

The elegant proof of~\eqref{Shao_thm} presented in~\cite{Sh97} basically reduced the problem to univariate LD problems for the projections of the vector $\bz :=(X, |X|^p)$ onto different directions $\bal=(\al_1, \al_2)\in \R_+\times \R_-$, by exploiting the observation that $x^{1/p}y^{1-1/p}
= p^{-1} \inf_{b>0}\Big(xb^{-1} + (p-1) y b^{1/(p-1)}\Big)$ for $x,y>0.$ That approach, however, neither provided a graphical clarification of why relation~\eqref{Shao_thm} actually holds nor explained  the ``true nature" of representation~\eqref{Shao_thm_1} for the limiting expression.

In the present paper, we directly relate the result stated in~\eqref{Shao_thm} to the standard  LD theory by reducing  the problem for self-normalized sums to a classical LD problem for an auxiliary bivariate random walk.
This clarifies the nature of the expression on the right-hand side of~\eqref{Shao_thm_1} for the limit and leads to an alternative proof of~\eqref{Shao_thm} under somewhat more general conditions. Moreover, this approach enables one to obtain exact (rather than just ``logarithmic") asymptotics for $\pr (W_{n,p}\ge z)$ as well. We also discuss possible extensions of these results to more general self-normalized settings.


The paper is structured as follows. Section~\ref{Sect_2} presents the above-mentioned reduction to the bivariate case, recalls necessary elements of the LD theory for multivariate random walks, and demonstrates an alternative proof of~\eqref{Shao_thm}.   Section~\ref{Sect_3} contains a derivation of the exact asymptotics for $\pr (W_{n,p}\ge z)$.  Section~\ref{Sect_4} discusses more general self-normalized settings.

\section{Reduction to the bivariate LD problem}
\label{Sect_2}

For $z>0,$ introduce the planar set
\begin{align}
\label{B_z}
B_z:= \{\bx=(x_1, x_2)\in \R^2: x_1\ge zx_2^{1/p}, \ x_2\ge 0\}
\end{align}
and consider the vector-valued random walk
\begin{align}
\label{Z_n}
\bZ_n:=\sum_j^n \bz_j, \quad n=1,2,\ldots,
\end{align}
with i.i.d.\ jumps
\begin{align}
\label{z_n}
\bz_j:= (X_j, |X_j|^p), \quad j=1,2,\ldots
\end{align}
The case $X\le 0$ being trivial, we will assume throughout this paper that $\pr (X>0)>0.$

Our key initial  observation is that, using the convention $0/0=\infty$ from~\eqref{Wnp}, the event in the probability on the left-hand side of~\eqref{Shao_thm} is equal to
\begin{align}
\label{redu_2}
A_n :=
 \biggl\{\frac{S_n}n \ge z \biggl(\frac{T_n}n\biggr) ^{1/p}\biggr\}
 = \biggl\{\bigg(\frac{S_n}n, \frac{T_n} n\biggr)\in B_z  \biggr\}
 = \biggl\{ \frac{\bZ_n}n \in B_z  \biggr\}.
\end{align}

This reduces the problem of evaluating the probabilities for the self-normalized sums on the left-hand side of~\eqref{Shao_thm} to the standard LD problem for the bivariate random walk $\{\bZ_n\}$. The latter problem one can attempt to solve using the by-now classical results in  the multivariate  setting. Assuming for the time being that $\bz_j\deq\bz $ are general i.i.d.\ random vectors  in $\R^d,$ $d\ge 2,$  we will now recall some relevant key  concepts and facts from the classical LD theory for random walks of the form~\eqref{Z_n} with such jumps that are essential for understanding the further exposition below.


Denote   by 
$\|\cdot \|:= \sqrt{\bx \bx^T }$ the Euclidean norm of $\bx\in \R^d$, $\bx^T$ being  the transpose of the (row-)vector~$\bx,$
and denote by $(B)$ and $[B]$ the interior and closure of the set $B\subset \R^d,$ respectively. We will now assume  that the jump distribution is non-degenerate in the sense that
\begin{align}
\label{non-deg}
\pr ( \bla \bz^T=c)<1 \quad \mbox{for any $\bla  \in \R^d $ and $c\in \R.$}
\end{align}
This condition is always met in the case of a random walk with jumps~\eqref{z_n}, except  when~$\bz$ has a two-point distribution\,---\,but then the problem is reduced to a straightforward univariate case to be discussed in detail in Section~\ref{Sect_3} below, in the context of deriving exact asymptotics of the LD probabilities.

Next let
 \begin{align}
\label{set_A}
\cA:=\{ \bla
 \in \R^d : A(\bla):= \ln \exn e^{ \bla \bz^T }
 <\infty \} .
\end{align}
It is well-known that  the function~$A$  and the set~$\cA$ are both convex (see e.g.~\cite{BoRo65}), and $A$ is analytic in~$(\cA)$. Clearly,  $\boldsymbol{0}\in \cA$.

For $\bla\in \cA,$ consider a random vector $\bz^{(\bla)}$ with the ``tilted distribution" (a.k.a.\ the Cram\'er transform, or the conjugate, of the distribution of~$\bz$ at the point~$\bla$) given by
 \begin{align}
\label{Cramer_transf}
\pr (\bz^{(\bla)}\in d\bx)
=  e^{ \bla \bx^T  -A(\bla)}
\pr (\bz \in d\bx).
\end{align}
One can easily see that, for any $\bla \in (\cA),$ $\bz^{(\bla)}$ has a finite exponential moment and that the mean of that vector is equal to $\exn\bz^{(\bla)}=\mbox{grad\,} A(\bla) .$ We set
\[
\cA':= \{\mbox{grad\,} A(\bla) : \bla \in (\cA)\}.
\]

Note that, in the special case when the jumps $\bz_j $ in our random walk are of the form~\eqref{z_n}, using the elementary bound
 \begin{align}
\label{xi<}
\la_1 x +\la_2 |x|^p\le \frac{p-1}p \frac{|\la_1|^{p/(p-1)}}{(p|\la_2|)^{1/(p-1)}}
\end{align}
valid for all $\la_2<0$ and $\la_1, x\in\R,$ one concludes that $\R\times (-\infty,0)\subseteq  (\cA),$ meaning that the following form of the Cram\'er condition   is met:
 \begin{align}
\label{Cramer_cond}
(\cA)\neq \varnothing .
\end{align}

Consider the rate function~$\La$ for $\bz$, defined as the Legendre transform (a.k.a.\ the conjugate) of the cumulant function~$A$:
\begin{align}
\label{Lambda}
\La (\bal):=\sup_{\bla \in \R^d} ( \bal  \bla^T  - A(\bla))\le \infty, \quad \bal \in \R^d.
\end{align}
The function~$\La\ge 0$ is convex and lower semi-continuous, and is analytic on~$\cA'.$ For $\bal \in \cA',$  the supremum on the right-had side of~\eqref{Lambda} is attained at the point
\begin{align}
\label{lambda}
\bla (\bal):= \mbox{grad\,}  \La (\bal)
\end{align}
(for these and further properties of the rate function, see e.g.\ Section~2 in~\cite{BoMo92}, Section~1.2 in~\cite{Bo20}, and also~\cite{BoMo78}, \cite{BoMo98}; for a detailed discussion of the  conjugates of  more general convex functions, see Section~12 in~\cite{Ro70}).

The probabilistic meaning of $\La$  is as follows (property~(5) in~\cite{BoMo92}): for any $\bal\in \R^d,$
 \begin{align}
\label{local_LD}
\La (\bal) = -\lim_{\ep \to 0} \lim_{n\to \infty }\frac1n \ln \pr \big(n^{-1}\bZ_n\in (\bal)_\ep \big),
\end{align}
where $(\bal)_\ep:=\{\bx \in \R^d: \|\bx - \bal  \|<\ep\},$ $\ep >0.$

For a set $B\subset \R^d$, put $\La(B):= \inf_{\bal\in B} \La (\bal).$ It immediately follows from~\eqref{local_LD} that, for any Borel set~$B$,
 \begin{align}
\label{lower_LD}
\liminf_{n\to \infty }\frac1n \ln \pr  (n^{-1}\bZ_n\in B)\ge -\La ((B)).
\end{align}
Moreover, for any {\em bounded\/} Borel~$B,$ one also has
 \begin{align}
\label{upper_LD}
\limsup_{n\to \infty }\frac1n \ln \pr  (n^{-1}\bZ_n\in B)\le -\La ([B])
\end{align}
(see e.g.\ p.~82 in~\cite{BoMo92}). Under the additional assumption
\begin{align}
\label{0_in_A}
\boldsymbol{0}\in (\cA),
\end{align}
the upper bound~\eqref{upper_LD} holds true for unbounded sets as well (see e.g.\ p.~84 in~\cite{BoMo92} or Corollary~6.1.6 in~\cite{DeZe98}). However,  generally speaking, this is not so  when~\eqref{0_in_A} is not met~(\cite{Di91,Mo14}), which is the case in general for our special  jumps~\eqref{z_n} (note that~\eqref{0_in_A} implies~\eqref{Cramer_cond}, but not the other way around).

Remarkably enough, it was shown in~\cite{Mo14} that, in the bivariate case $d=2$, relation~\eqref{upper_LD} always  holds for unbounded Borel sets~$B$ as well. In the case $d=3,$ however, there exist a random vector~$\bz$ satisfying~\eqref{Cramer_cond} (but not~\eqref{0_in_A}) and an unbounded set~$B$ such  that~\eqref{upper_LD} fails (Theorem~1.4 in~\cite{Mo14}).

It follows that, in the general bivariate case, for any Borel set~$B\subset  \R^2,$ both~\eqref{lower_LD} and~\eqref{upper_LD} hold true, so that the sequence of the distributions of $n^{-1}\bZ_n,$ $n\ge 1,$ satisfies the so-called LD principle with the rate function~$\La$. One consequence of this is that if $\La ((B))=\La([B]) $  then  there always exists the limit
 \[
 \lim_{n\to \infty }\frac1n \ln \pr  (n^{-1}\bZ_n\in B) =  -\La (B).
 \]

Now we  will turn back to the special case of random vectors~\eqref{z_n} and events~\eqref{redu_2}. Note that the set~$B_z$ from~\eqref{B_z} is clearly closed and unbounded. However, it is not hard to see that, in the case when $\pr (X=0)>0,$ one can have $\La ((B_z))>\La([B_z])$. For instance,   $\La ((B_1))=\infty$ (due to~\eqref{S<T}) but $\La([B_z])\le \La (\boldsymbol{0}) =-\ln \pr (X=0)$, $z>0.$ 

Hence, in the case of such sets, the direct application of the LD principle by just computing $\La (B_z)=\La ([B_z])$ (note that $B_z$ is closed) is not feasible. One actually needs to ``tighten" the bounds (as we will see, this applies to the lower one only). This makes the problem somewhat more challenging.

Let $\widetilde{B}_z:= \{(x_1, x_2)\in \R^2: x_1 =  zx_2^{1/p}, \ x_2\ge 0\} $  be the left-upper component of~$\partial B_z$. The following result  shows that the event in~\eqref{redu_2} is indeed an LD for the sum~$\bZ_n$ provided that $z>z^*$ (see~\eqref{z*}).

\begin{lemo}
\label{prop}
For $z\in (z^*,1)$, one has
$ \La ( B_z )=\La (\widetilde{B}_z)<\infty.$
\end{lemo}

\begin{rema}
\label{rem_2}
\rm
As we already noted, $\La ((B_1))=\infty$, so that one clearly still has $ \La ( B_z )=\La (\widetilde{B}_z)$ for $z\ge 1$ as well, the common value being $\La (\boldsymbol 0)= -\ln \pr (X=0)\le \infty$ for $z>1$ (cf.~Remark~\ref{rem_1}), whereas for $z=1$ it may be different when the distribution of~$X$ has atoms in~$(0,\infty).$
\end{rema}

\noindent
{\em Proof of Lemma~\ref{prop}.} First we will show that $\cA'\cap  {B}_z \neq \varnothing.$ Since, as we already mentioned,  $\La (\bal )<\infty$ (and, moreover, is analytic) at any~$\bal\in  \cA'$, this will prove that $\La ({B}_z)<\infty.$

\begin{figure}[ht]
\begin{center}
	\includegraphics[scale=.4]{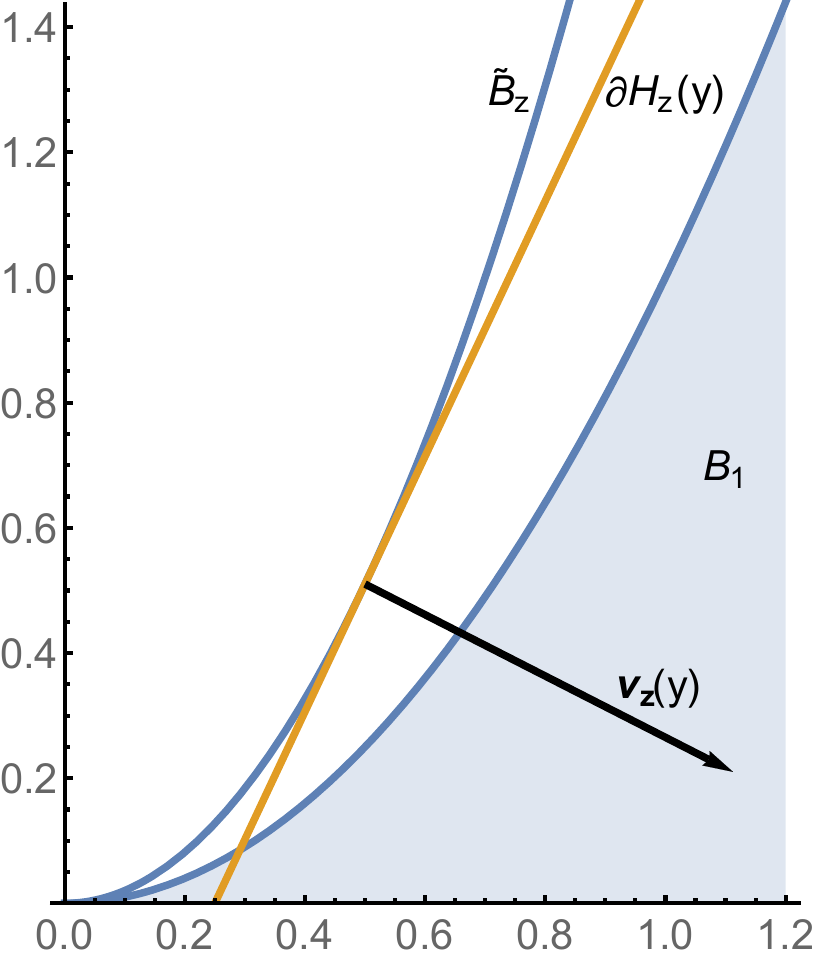}
\end{center}
	\caption{An illustration to the definitions of $B_1,$ $\widetilde{B}_z,$ $\bnu_z (y),$ and $H_z (y)$. The vector $\bnu_z (y)$ is shown not to scale (it is much longer for the chosen values $p=2,$ $z=0.7,$ $y=0.5$.)}
	\label{Fig_1}
\end{figure}

For $z>0,$ $y\ge 0,$ introduce the vector
%
\[
\bnu_z (y):=  (py^{p-1}z^{-p},-1).
\]
It is easily seen that $\bnu_z (y)$ is  orthogonal to the tangent to $\widetilde{B}_z$ at the point $\by_z : = (y,y^pz^{-p})\in \widetilde{B}_z$ and points inside~${B}_z$.  Set
\[
h_z (y):=  \by_z \bnu_z^T  (y) =  (p-1) y^pz^{-p}
\]
and introduce  the closed half-planes
\begin{align}
\label{Hzy}
H_z (y):= \{\bx=(x_1, x_2) :   \bx \bnu_z^T (y)
 \ge h_z (y)\}
\end{align}
(see Fig.~\ref{Fig_1}).

Since  $\pr (X>0)>0,$   there exists a $y>0$ such that $\pr (X\in [y, y+\delta])>0$ for any~$\delta >0.$ As the support of the distribution of~$\bz$ is part of the convex set
\[
 B_1^\cup :=\{\bx: |x_1| \le x_2^{1/p}, x_2 \ge 0\}\subset [H_{y,1}^c],
\]
$\widetilde{B}_1$ being the ``right branch" of~$\partial  B_1^\cup,$    and $  (y, y^p) \bnu_1^T  (y)    = (p-1)y^p= h_1 (y) , $ we conclude that, for
the random variable $\xi:=  \bz \bnu_1^T  (y)  ,$ one has
\[
\pr (\xi \le h_1 (y))=\pr (\bz  \in [H_{y,1}^c] ) =1, \quad
p_\delta := \pr (\xi \in [h_1 (y) - \delta, h_1 (y) ])>0
\]
for any $\delta >0.$

Denote by $\xi^{(\la)},$ $\la \ge 0,$ the  Cram\'er transform  of $\xi$, i.e., a random variable with density $e^{\la x}/\exn e^{\la \xi}$ w.r.t.\ the distribution of~$\xi.$ One has
\[
\pr (\xi^{(\la)} <h_1 (y) - 2 \delta )
 = \frac{\exn ( e^{\la \xi}; \xi <h_1 (y) - 2 \delta)}{\exn e^{\la \xi}}
 \le  \frac{  e^{\la  (h_1 (y) - 2 \delta)} }{p_\delta  e^{\la  (h_1 (y) -   \delta)}}
 \to 0 , \quad \la\to\infty.
\]
Since the family $\{\xi^{(\la)}\}_{\la \ge 1}$ is evidently uniformly integrable (cf.~\eqref{xi<}), we conclude that $\exn \xi^{(\la)}\to h_1 (y)$ as $\la\to\infty.$ That is, for any $\ep >0,$ there is a $\la_\ep <\infty $ such that
\begin{align}
\label{Exi}
\exn \xi^{(\la)}>h_1 (y)-\ep \|\bnu_1 (y )\| \quad \mbox{for $\la >\la_\ep.$}
\end{align}

Now clearly $\xi^{(\la)} \deq   \bz^{(\la \bnu_1 (y))}  \bnu_1^T  (y)     $  and $\la \bnu_1 (y) \in (\cA),$ $\la >0.$
Hence, letting
\begin{align}
\label{bold_e}
\be_z (y) : = \bnu_z (y) \|\bnu_z (y )\|^{-1}, \quad y\ge 0, \ z\in (0,1],
\end{align}
be the unit vector collinear with $\bnu_z (y)$ and noting that $\exn \bz^{(\la \bnu_1 (y))}\bnu_1^T (y) = \exn \xi^{(\la)},$ we obtain from~\eqref{Exi} and~\eqref{Hzy} that
\[
\exn \bz^{(\la \bnu_1 (y))} \in H_{y,1} -\ep \be_1  (y)  ,
\quad \la > \la_\ep.
\]
On the other hand, since  $\pr (\bz^{(\la \bnu_1 (y))}\in B_1^\cup)= \pr (\bz \in B_1^\cup)=1 $, where $B_1^\cup $ is strictly convex, and the distribution of~$\bz$ is non-degenerate, one has
\[
\exn \bz^{(\la \bnu_1 (y))} \in (B_1^\cup), \quad  \la >0.
\]
Since $\exn \bz^{(\la \bnu_1 (y))} \in \cA'$ and, for any $z\in (0,1),$
\[
(H_{y,1} -\ep \be_1  (y))\cap (B_1^\cup) \subset B_z
\]
for all small enough~$\ep>0,$ the desired assertion that $\cA'\cap B_z \neq \varnothing$ is proved.

To prove the claim $\La (B_z) = \La (\widetilde B_z),$ we will consider two alternative cases.

{\bf Case~1:}  $\exn |X|^p<\infty.$ In this case,  the random vector~\eqref{z_n} has finite mean $\boldsymbol{m}:=\exn \bz,$ the condition $z>z^*$ implying  that \begin{align}
\label{m_in_B}
\boldsymbol{m} \in (B_z^{\cup}),
\qquad B_z^{\cup}:=\{\bx\in \R^2: |x_1|\le zx_2^{1/p}, x_2 \ge 0\}.
\end{align}

It follows  from the definition of the rate function that $\La (\boldsymbol m)=0$. Hence $\La\ge 0$ attains its global minimum at the point~$\boldsymbol{m}$. Being convex, $\La$~is non-decreasing along any of the straight line  rays emanating from~$\boldsymbol{m}$. Any such ray that ``hits" $ {B}_z $ clearly enters that set at a point from~$\widetilde{B}_z,$ which implies that  $ \La ( B_z )=\La (\widetilde{B}_z).$

{\bf Case~2:}  $\exn |X|^p = \infty.$ To deal with this case, introduce truncated random variables $X^{[n]}:= X\ind_{\{|X|<n\}}$ and vectors $\bz^{[n]}:= (X^{[n]}, |X^{[n]}|^p),$ $n=1,2,\ldots$ These vectors clearly have finite means $\boldsymbol{m}^{[n]}=({m_1}^{[n]},{m_2}^{[n]}):=\exn \bz^{[n]}.$ If we knew that ${m_1}^{[n]}=o\big(({m_2}^{[n]})^{1/p}\big)$ as $n\to \infty,$ it would mean that these mean vectors will eventually lie in~$B_z^{\cup}.$ The  next lemma, of which the claim may be known, implies that this is so indeed.

\begin{lemo}
\label{lemo_1}
Let $p>1$ and $F$ be a distribution function on $[0,\infty)$ such that $\int_0^\infty t^p dF(t)=\infty.$ Then $\big( \int_0^x t  dF(t)\big)^p
=  o\big( \int_0^x t^p  dF(t)\big)$ as $x\to \infty.$
\end{lemo}

\noindent
{\em Proof.} Set $G(x):= \int_0^x t  dF(t).$ If $G(\infty)<\infty$ then the claim of the lemma is trivial. Assume that $G(x)\to\infty$ as $x\to \infty.$

Setting $\overline{F}(t):=1- F(t),$ we first note that, for any fixed $v>0,$
\begin{align}
\label{G_bound}
G(vx) = \int_0^{\sqrt{vx}} + \int_{\sqrt{vx}}^{vx}
\le \sqrt{vx} + vx  \overline{F}(\sqrt{vx}), \quad  xv\ge 1.
\end{align}
Clearly $v  \overline{F}(\sqrt{vx})\to 0$ as $x\to \infty.$ Hence there exists a  function $v_x\uparrow \infty,$ $v_x=o(x)$ as $x\to\infty$, that    increases slowly enough  to ensure that   $v_x  \overline{F}(\sqrt{v_x x})\to 0$. It follows from~\eqref{G_bound} that $G(v_x x) = o(x) $ and hence also $G(v_{G(x)} G(x) ) =o(G(x))$. Since clearly $G(x)=o(x)$, we can assume without loss of generality that $v_{G(x)} G(x)=o(x)$. Now  we have
\begin{align*}
\int_0^x t^p  dF(t)
 & =
 \int_0^x t^{p-1}  dG(t) \ge \int_{v_{G(x)} G(x)}^x t^{p-1}  dG(t)
 \\
 & \ge  (v_{G(x)} G(x))^{p-1}\big[G(x) - G(v_{G(x)} G(x) ) \big]
 \\
 & = v^{p-1}_{G(x)} G(x)^p (1+o(1)) \gg G(x)^p
\end{align*}
as $x\to \infty.$ Lemma~\ref{lemo_1} is proved. \hfill$\Box$

\medskip

Denoting by $\La^{[n]}$ the rate functions for $\bz^{[n]},$ it follows from Case~1 above that $ \La^{[n]} ( B_z )=\La^{[n]} (\widetilde{B}_z)$ for all large enough~$n$. Theorem~1 on p.~73 in~\cite{BoMo92} states that, under conditions~\eqref{non-deg} and~\eqref{Cramer_cond}, $ \La^{[n]} (\bal) \to \La  (\bal)$ at any continuity point~$\bal$ of~$\La$. Therefore the desired property $ \La  ( B_z )=\La  (\widetilde{B}_z)$ will be ``inherited" by the limiting rate function~$\La$ as well.

Lemma~\ref{prop} is proved.
\hfill$\Box$

\medskip

Now we will turn to  computing  the supremum and infimum  limits  for $n^{-1}\ln \pr (A_n),$ $n\ge 1,$ for events $A_n $ from~\eqref{redu_2}. As the set~$B_z$ is closed, in view of~\eqref{upper_LD} one has
\begin{align}
\label{La_H}
\limsup_{n\to\infty} n^{-1}\ln \pr (A_n )
\le - \La (B_z)
 = -  \inf_{y\ge 0 } \La (H_z (y)),
\end{align}
where the last equality follows from the observations that $B_z\cup (H_{z}(0)\setminus \{\boldsymbol{0}\})  =\bigcup_{y\ge 0} H_z (y)$ and $\La ( H_{z}(0)\setminus \{\boldsymbol{0}\})= \infty$ due to $\La (\bal) =\infty$ for all $\bal \in H_{z}(0)\setminus \{\boldsymbol{0}\}$ in view of~\eqref{local_LD}.

Using the same argument as we employed  in the proof of Lemma~\ref{prop} to show that $\La (B_z) = \La (\widetilde B_z)$, one obtains that $\La (H_z (y)) =  \La (\partial H_z (y))$. By the rate function property~(13) on p.~66 in~\cite{BoMo92}, one has
\[
\La (\partial H_z (y))\equiv \La (\{\bx: \bx  \bnu_z^T  (y)= h_z (y)\})=\La_{\xi_z(y)} (h_z (y)),
\]
where
\[
\xi_z (y):=   \bz \bnu_z^T  (y) =  py^{p-1} z^{-p} X - |X|^p
\]
and
\begin{align*}
\La_{\xi_z(y)}  (a) :  = \sup_{\la\in \R}(a \la   -\ln \exn e^{\la \xi_z(y)})
  = - \inf_{\la\in \R} \ln \exn e^{\la ( \xi_z(y) - a )}, \quad a\in \R.
\end{align*}
We obtain that
\begin{align}
\La (H_{ y ,z})  &  =\La_{\xi_z(y)} (h_z (y))
\notag
\\
\label{La_xi}
&
 = - \inf_{\la\in \R} \ln \exn e^{\la ( \xi_z(y) - h_z (y))}
 = - \inf_{\la\ge 0} \ln \exn e^{\la ( \xi_z(y) - h_z (y))},
\end{align}
where we   replaced $\inf_{\la\in \R}$ with $\inf_{\la\ge 0}$ since the infimum is attained on the positive half-line due to the   inequality
\begin{align}
\label{<0}
\exn \xi_z(y)  <h_z (y).
\end{align}
To verify~\eqref{<0}, we observe  that  if $\exn |X|^p=\infty$ then $\exn \xi_z(y)  = -\infty$, while if $\exn |X|^p<\infty $  then one has
\[
\exn (\xi_z (y)- h_z(y))=\exn     \bz \bnu_z^T  (y) - h_z(y)
 =       \boldsymbol{m} \bnu_z^T  (y)  - h_z(y) <0
\]
due to~\eqref{m_in_B} and the obvious inclusion~$(B_z^\cup)\in H_z (y)^c =\{\bx: \bx \bnu_z^T  (y)   <h_z(y)\}.$ Together with~\eqref{La_H} this yields
\begin{align}
\limsup_{n\to\infty} n^{-1} & \ln \pr (A_n )
\notag
\\
& \le
\sup_{y\ge 0 }\inf_{\la\ge 0} \ln \exn \exp\{\la [ py^{p-1} z ^{-p} X - |X|^p-(p-1) y^pz ^{-p} ] \}.
\label{P=La}
\end{align}
Letting here $\la:=zt/p$ and $y:=z c^{1/(p-1)},$ we obtain the quantity~\eqref{Shao_thm_1} appearing on right-hand side of~\eqref{Shao_thm}.

To get the lower bound, note that clearly $B_z\supset H_z (y)\cap \{x_2\ge 0\}$ for any $y\ge 0,$ implying
\begin{align*}
\pr (n^{-1} \bZ_n \in B_z)
 & \ge \sup_{y\ge 0}\pr (n^{-1} \bZ_n \in  H_z (y))
 \\
 & = \sup_{y\ge 0}\pr \Bigl(\frac1n \sum_{j=1}^n   \xi_{z,j} (y)
  \ge  h_z (y)\Bigr),
\end{align*}
where $\xi_{z,j} (y): =  \bz_j  \bnu_z^T  (y) $ are i.i.d.\ copies of~$\xi_{z} (y).$ By the ``strengthened version" of Cram\'er's theorem (see  Corollary~2.2.19 in~\cite{DeZe98}),  there exists the limit
\begin{align*}
\lim_{n\to\infty} \frac1n \ln \pr   \Bigl(\frac1n\sum_{j=1}^n   \xi_{z,j} (y) \ge  h_z (y)\Bigr)
    = - \inf_{v\ge h_z (y)} \La_{\xi_z (y) } (v)
  = -   \La_{\xi_z (y) } (h_z (y)),
\end{align*}
where the last equality is due to~\eqref{<0}.
Therefore, using~\eqref{La_xi} and~\eqref{<0}, we get
\begin{align*}
\liminf_{n\to\infty} \frac1n   \ln \pr (A_n )
&
\ge
\liminf_{n\to\infty} \sup_{y\ge 0} \frac1n \ln \pr \Bigl(\frac1n \sum_{j=1}^n   \xi_{z,j} (y)  \ge  h_z (y)\Bigr)
\\
&
\ge \sup_{y\ge 0}(-   \La_{\xi_z (y) } (h_z (y)))
=  \sup_{y\ge 0} \inf_{\la \ge 0} \ln \exn e^{\la ( \xi_z(y) - h_z (y))},
\end{align*}
 which coincides with the right-hand side of~\eqref{P=La}.

Thus, directly based on the bivariate LD theory, we proved the following result extending Theorem~\ref{thm_Shao}.

\begin{theo}
\label{thm_our_1}
The claim of Theorem~{\rm \ref{thm_Shao}} holds in the general case for   $z>z^*,$  without  the assumption that   either  $\exn X\ge 0$ or $\exn |X |^p =\infty.$
\end{theo}

\begin{rema}
{\em
It was noted in Remark~1.2 in~\cite{Sh97} that~\eqref{Shao_thm} remains true in the special case $p=2$, $\exn X^2<\infty$  without the assumption~$\exn X\ge 0.$ We showed that the latter assumption is irrelevant  in the general case $p>1$ as well.  }
\end{rema}

\section{Exact asymptotics}
\label{Sect_3}

The LD probability approximation $\pr(A_n)=e^{ n J_z + o(n)}$ established in~\eqref{Shao_thm} is, of course, very crude due to the multiplicative error term~$e^{o(n)}.$ Our reduction~\eqref{redu_2} of the  LD problem for self-normalized random walks to the classical LD problem for bivariate random walks with jumps~\eqref{z_n} enables one to obtain much sharper results. We will present them in this section, the form of the approximation depending on whether or not the jumps satisfy the non-degeneracy condition~\eqref{non-deg}.

\medskip
\noindent
{\bf Case 1.} First consider the case where the non-degeneracy condition~\eqref{non-deg} is not met. This can only happen when there exist $a<b$ such that
\begin{align}
\label{2_point_distr}
1-\pr (X=a) = \pr (X=b) =:q \in (0,1).
\end{align}
The case $b\le 0$ is trivial (note that $\pr (A_n)=0$ when $b<0$ and $\pr (A_n)=\pr (T_n=0)= q^n$ when $b=0$), so we will assume that $b>0$.

Set
\[
f(t) := \frac{a + t (b-a)}{(|a|^p + t (b^p-|a|^p))^{1/p}}, \quad t\in [0,1], \]
and observe that $f(q) =z^* = \exn X(\exn |X|^p)^{-1/p}.$

Under the assumption that $z>z^*,$ if  $a<0$ then there exists a unique solution $t_z\in (0,1)$ to the equation $f(t)= z$, while if $a\ge 0$ then there exist two solutions $0\le t_{z-} <t_{z_+}<1$ to that equation, where $t_{z-}=0$ iff~$a=0.$  In both cases, the respective points of the form
\[
 (a + t_z (b-a), |a|^p + t_z (b^p-|a|^p))\in \R^2
\]
with $t=t_z$ in the former case and $t=t_{z\pm}$ in the latter one, are just the points where the straight line segment~$I$ connecting the points $(a, |a|^p)$ and $(b, b^p)$  crosses the curve~$\widetilde{B}_z$.

Setting $N_n:=\sum_{j=1}^n \eta_j,$ where $\eta_j:=\ind_{\{X_j=b\}},$  $j=1,2,\ldots,$ are i.i.d.\ Bernoulli random variables with success probability~$q$, we see that $S_n= na + N_n (b-a),$ $T_n= n|a|^p + N_n (b^p-|a|^p),$ so that one always has $n^{-1}\bZ_n\in I$, $ n\ge 1,$ and, moreover,
\begin{align}
\label{A_n_N_n}
A_n= \left\{
\begin{array}{ll}
\{N_n \ge t_z n\}, & \quad a <0,
\\
\{N_n = 0\}\cup \{N_n\ge t_{z+} n\}, & \quad a = 0,
\\
\{N_n\le t_{z-} n\}\cup \{N_n\ge t_{z+} n\}, & \quad a > 0.
\end{array}
\right.
\end{align}
The assumption that $z>z^*=f(q)$ implies that $t_z>q=\exn \eta $ when $a<0$, and that $t_{z-}<q <t_{z+}$ when $a \ge 0.$   That is, the event $A_n$ is equivalent to the respective combinations of LD events for the binomial random variable~$N_n.$

To evaluate the probabilities on the events on the right-hand side of~\eqref{A_n_N_n}, we have to find the rate function  $\La_\eta$ for $\eta:=\ind_{\{X=b\}}.$
Direct computation yields that, for $ \al \in (0,1)$,
\[
\La_\eta (\al ) = \sup_{\la \in \R}(\al \la - \ln (1-q + qe^\la)) = \al\ln\frac{\al}{q}+(1-\al)\ln\frac{1-\al}{1-q},
\]
where the supremum is attained at the point $\la (\al):= \ln\frac{\al}{1-\al}\frac{1-q}q.$ Note also that, as one could expect from~\eqref{local_LD}, one has   $\La_\eta  (0)=-\ln (1-q),$ $\La_\eta (1) =-\ln q$, these values obtained as $\la\to \mp\infty,$ respectively.

For the Cram\'er transform $\eta^{(\la  )}$ of $\eta$ with the above value  $\la = \la(\al)$, one  has
\begin{align*}
\pr (\eta^{(\la (\al))}=0)
& =\frac{1-q}{\exn e^{\la (\al )\eta}}
=\frac{1-q}{1-q + qe^{\la (\al )}}
\\
& =\frac{1-q}{1-q + (1-q)\al /(1-\al)} = 1-\al,
\end{align*}
so that $\pr (\eta^{(\la (\al))}=1) =\al $ and hence $\mbox{Var} (\eta^{(\la (\al))}) = \al (1-\al). $

As $N_n$ is integer-valued, for stating the exact LD probability asymptotics, we will need the following auxiliary quantities:
\[
\al_z:=\lceil t_z\rceil ,
\quad \al_{z-}:=\lfloor  t_{z-}\rfloor ,
\quad \al_{z+}:=\lceil t_{z+}\rceil ,
\]
where $\lfloor x\rfloor:=\max\{ k\in \Z: k\le x\} $ and  $\lceil x\rceil:=\min\{ k\in \Z: k\ge x\} $ for $x\in \R.$

Now from~\eqref{A_n_N_n} and Corollary 6.1.7 in~\cite{BoBo08} we obtain the following exact asymptotics results.
\begin{theo}
 \label{thm_Exact_degen}
Assume that $p>1$ and  $X$ has the two-point distribution~\eqref{2_point_distr} with $b>0$.
Then, for $z\in (f(q),1),$ as $n\to\infty,$ one has$:$
\smallskip

{\rm (i)}~if $a< 0$ then
\begin{align*}
\pr ( W_{n,p}
   \ge z  )
    = \frac{(1+o(1))e^{-n \La_\eta (\al_z)}}{(1- e^{-\la (\al_z)})\sqrt{2\pi n \mbox{\rm Var} (\eta^{(\la (\al_z))})}}
    = (1+o(1)) Q_n (\al_z),
\end{align*}
where
\[
Q_n(\al):= \sqrt{\frac{\al }{ 2\pi n( 1-\al )} }\frac{1-q}{\al  -q}
   \Bigl(\frac{q}{\al }\Bigr)^{n\al } \Bigl(\frac{1-q}{1-\al }\Bigr)^{n(1-\al )}, \quad \al \in (0,1) ;
\]

\smallskip

{\rm (ii)}~if $a=  0$ then
\begin{align*}
\pr ( W_{n,p}
   \ge z  )
    = (1-q)^n +  (1+o(1)) Q_n (\al_{z+});
\end{align*}

\smallskip

{\rm (iii)}~if $a> 0$ then
\begin{align*}
\pr ( W_{n,p}
   \ge z  )
    = (1+o(1)) (Q_n (\al_{z-}) +  Q_n (\al_{z+})).
\end{align*}
\end{theo}

\medskip
\noindent
{\bf Case 2.} Now we assume that condition~\eqref{non-deg} is met. In this case, the desired asymptotics can   be obtained from Theorem~2 on p.~98 in~\cite{BoMo92}.   The conditions ensuring the validity of the claim of that theorem are listed in  Section~2 of \S\,4 of that paper. One of them is~\eqref{non-deg}.

Fix an arbitrary $z\in (z^*, 1).$ The next assumption of the above-mentioned theorem from~\cite{BoMo92} is that there exists a {\em unique\/} point
\[
\wbal= ( \widehat{\al}_1, \widehat{\al}_2)=\wbal_z\in B_z
\]
such that $\La (\wbal) = \La (B_z)$ and $\widehat{\al}_1>0$. It follows from Lemma~\ref{prop}  that      $ \wbal\in \widetilde B_z$.

\begin{figure}[ht]
\begin{center}
	\includegraphics[scale=.55]{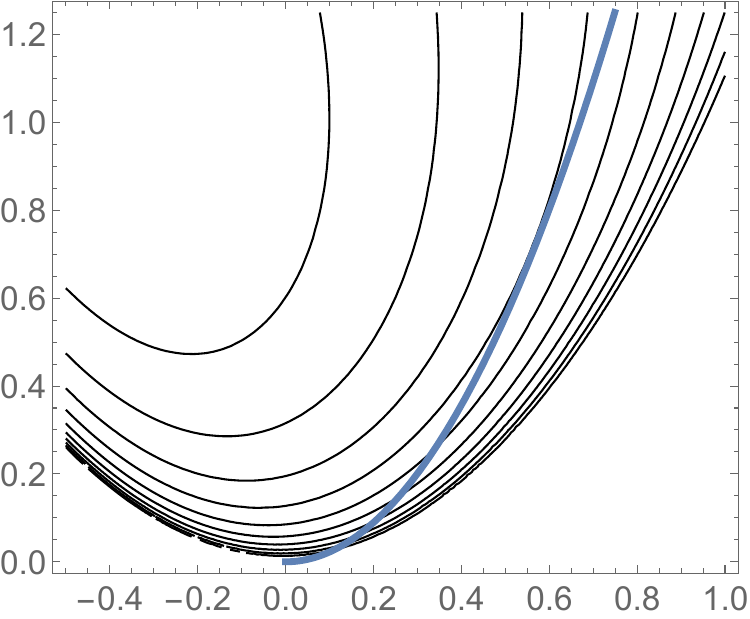}
\end{center}
	\caption{The thick line is $\widetilde{B}_z$ with $p=2,$ $z=0.67.$ The contour plot is that of the rate function $\La = \La_{\bz}$ for $\bz =(X, X^2)$ with $X\sim N(-0.5,1).$ The ``first order only" contact of~$\widetilde{B}_z$ is with the level line~$\La (\bal)=\ell$ with $\ell \approx 0.72.$}
	\label{Fig_2}
\end{figure}

Clearly, the level line  $L_z:=\{\bal\in \R^2: \La (\bal)= \La (\wbal)\}$ must be ``tangent" to the boundary $\widetilde{B}_z$ at the point~$\wbal$. The next condition we need is stated as~(6) on p.\,94 in~\cite{BoMo92} (equivalently, as~(11) or~(12) 
on p.\,95 therein). It states that the ``contact" between the lines~$L_z$ and $\widetilde{B}_z$ at~$\wbal$ must be of the first order only, meaning that the quadratic approximation to~$L_z$ in vicinity of~$\wbal$ should differ from the one for the curve~$\widetilde{B}_z$. Fig.~\ref{Fig_2} shows an example of this type of contact.

Note that $\widetilde{B}_z$ is given by $\al_2 =g_1(\al_1):= z^{-p}\al_1^p,$ $\al_1\ge 0,$  and, by the implicit function theorem, there exists a $\de>0$ and a smooth function~$g_2$ such that, in vicinity of~$\wbal,$ the line~$L_z$ is given by  $\al_2 =g_2(\al_1),$ $\al_1\in (\widehat{\al}_1 - \de, \widehat{\al}_1  + \de).$ By the definition of~$\wbal,$ one has $g_2(\al_1)\ge g_1(\al_1)$ in that neighborhood of $\widehat{\al}_1. $ Therefore the above-mentioned condition~(6) from~\cite{BoMo92} can be stated as
\begin{align}
\label{cond_C3}
g_2'' (\widehat{\al}_1) > g_1'' (\widehat{\al}_1).
\end{align}
This condition can be paraphrased in terms of the rate function~$\La.$ Indeed, that the lines~$\widetilde{B}_z$ and~$L_z$ have a common tangent at~$\wbal$ means that
\begin{align}
\label{g''}
g_2' (\widehat{\al}_1)=g_1' (\widehat{\al}_1)= p  z^{-p}\widehat{\al}_1^{p-1}. 
\end{align}
Differentiating twice the identity $\La (x, g_2(x))= \La (\wbal), $ $x\in (\widehat{\al}_1 - \de, \widehat{\al}_1  + \de),$ yields
\begin{align*}
 &   \partial_{\al_1 \al_1}  \La (x, g_2(x))
 + 2 \partial_{\al_1 \al_2}  \La (x, g_2(x)) g'_2(x)
 \\
  & \qquad\qquad\qquad + \partial_{\al_2 \al_2}  \La (x, g_2(x)) (g'_2(x))^2
 + \partial_{\al_2} \La (x, g_2(x)) g''_2(x)=0.
\end{align*}
Using~\eqref{g''}, we can now restate condition~\eqref{cond_C3} as
\begin{multline}
 -\frac1{\partial_{\al_2 } \La (\wbal)}
 \big(
 \partial_{\al_1 \al_1} \La (\wbal)  + 2\partial_{\al_1 \al_2} \La (\wbal)
 p  z^{-p}\widehat{\al}_1^{p-1} + \partial_{\al_2 \al_2} \La (\wbal)
( p  z^{-p}\widehat{\al}_1^{p-1})^2
 \big)
\\
 >
 p(p-1) z^{-p}\widehat{\al}_1^{p-2}  \quad \mbox{for $\wbal=\wbal_z$}.
 \label{Cond(6)}
\end{multline}

To state our main assertion in Case~2, it remains to compute the quantity $\chi (\cdots)$ appearing in the formulation of Theorem~2 on p.\,98 in~\cite{BoMo92}. To this end, denote by
\[
\Sigma(\bal):=\exn (\bz^{(\bla(\bal))}- \exn\bz^{(\bla(\bal))})^T (\bz^{(\bla(\bal))}-\exn\bz^{(\bla(\bal))} )
\]
the covariance matrix of the Cram\'er transform $\bz^{(\bla)}$ with $\bla = \bla (\bal)$ (see~\eqref{Cramer_transf}, \eqref{lambda}), let $\be (\bal):=\be_z ( {\al}_1) $ for $\bal=(\al_1,\al_2)\in \widetilde{B}_z$  (see~\eqref{bold_e}) and  denote by $\overline \be (\bal)$   a unit vector orthogonal to~$\be(\bal)$ (it does not matter which of the two possible directions is chosen). Using the observation that $\be^T ( \bal)\be(\bal)  + \overline \be^T  (\bal)  \overline \be (\bal)$ is the  $2\times 2$ identity  matrix (just check how this matrix acts as an operator), we find that  the covariance matrix~(9) from p.\,95 in~\cite{BoMo92} takes the form
\[
\widetilde \Sigma (\bal) := \overline \be^T  (\bal) \overline \be (\bal)
\bigg(
\Sigma(\bal)
-\frac{ \Sigma(\bal)\be^T (\bal) \be(\bal)\Sigma(\bal)}{\be(\bal)\Sigma(\bal)\be^T (\bal) }
\bigg)
\overline \be^T  (\bal) \overline \be (\bal).
\]
In the bivariate case, $\widetilde \Sigma (\bal)=(\widetilde\Sigma_{ij} (\bal))_{i,j=1,2}$ is the covariance matrix of the zero mean normal distribution concentrated on the straight line   $\{\bx \in \R^2: \bx \be^T (\bal) =0\}$ that was described in Lemma~10 on p.\,116 in~\cite{BoMo92}. Denoting by $\boldsymbol{\eta} (\bal)$ a random vector with that distribution, the desired quantity $\chi (\cdots)$ is defined as $\chi^* (\wbal),$ where
\[
\chi^* (\bal):=\exn \exp\bigg\{\frac{\|\bla (\bal)\|}2
 \boldsymbol{\eta}  (\bal)D (\bal) \boldsymbol{\eta}^T  (\bal)\bigg\} 
\]
and $D(\bal)$ is the Hessian matrix of the function $V(\bal):=z^{-p}\al_1^p -\al_2$ that specifies the curve~$\widetilde{B}_z$  as $\{\bal\in \R^2: V(\bal)=0, \al_1\ge 0\}$. Clearly, the only non-zero entry in this matrix is~$D_{11} (\bal)=\partial_{\al_1\al_1}V(\bal)=p(p-1)z^{-p}\al_1^{p-2}$. Therefore, denoting by~$Y$ a standard normal variable and setting
\[
\sigma^2(\bal) :=\|\bla (\bal)\|D_{11} (\bal) \widetilde\Sigma_{11} (\bal) ,
\]
we conclude that
\[
\chi^* (\bal)=\exn e^{\sigma^2(\bal)Y^2/2}=\frac1{\sqrt{1- \sigma^2(\bal)}},
\]
provided that $\sigma^2(\bal)<1.$ Note that condition~\eqref{Cond(6)} ensures that $\sigma^2(\wbal) <1$, see p.\,95 in~\cite{BoMo92}.

Recalling that $\La (\bal_z)=-J_z$ (see~\eqref{Shao_thm}, \eqref{Shao_thm_1} and our proof of Theorem~\ref{thm_our_1}), we can now state the following assertion as a direct consequence of Theorem~2 on p.\,98 in~\cite{BoMo92}.

\begin{theo}
\label{thm_Exact_nondeg}
Under the conditions stated for Case~$2$ in this section, for $z\in (z^*,1)$ one has
\[
\pr ( W_{n,p}    \ge z  ) =  \frac{(1+o(1)) e^{nJ_z}}%
{\sqrt{2 \pi n (1-   \sigma^2(\wbal_z))\be (\wbal_z)\Sigma (\wbal_z)\be^T (\wbal_z) }\|\bla (\wbal_z)\| }, \quad  n\to\infty.
\]
\end{theo}

\section{On more general self-normalized settings}
\label{Sect_4}

In conclusion, we will observe that our approach from Section~\ref{Sect_2} enables one to treat the LD problem for self-normalized random walks in  more general settings as well.

Namely, assume that $u:\R\to [0,\infty)$ is a convex function such that $u(0)=0 $ and both functions $u(\pm x),$ $x\ge 0,$ are strictly increasing, and let $u^{-1}_+(y),$ $y\ge0,$ be the inverse of $u(x),$ $x\ge 0.$ Next set $U_n:= \sum_{j=1}^n u (X_j),$ $n\ge 1,$ and consider   self-normalized sums of the form
\begin{align*}
W_{n}:=  \frac{S_n}{n u_+^{-1}(U_n/n)},
\end{align*}
where, similarly to~\eqref{Wnp}, $W_n:=\infty$ in the case when   $U_n=0.$

In the original problem dealt with in Theorem~\ref{thm_Shao}, one had $u(x)=|x|^p,$ $p>1,$ yielding $u_+^{-1}(y)=y^{1/p}$ and turning the above~$W_n$ into~$W_{n,p}$ from~\eqref{Wnp}.

For $z>0,$ introduce the planar set
\[
C_z: = \{\bx \in \R^2: x_1\ge zu_+^{-1}(x_2), \ x_2\ge 0\}.
\]
As in~\eqref{redu_2}, one clearly has
\[
\{W_n \ge z\} = \{(S_n/n, U_n/n)\in C_z\}=\{\bZ_n/n\in C_z\}
\]
for sums~\eqref{Z_n} with $\bz_j:= (X_j, u(X_j)),$ $j\ge 1.$ It is not hard to see that the above event will be an LD when
\[
z> z^{**}:=\left\{
\begin{array}{ll}
\exn X/u_+^{-1} (\exn u(X))  & \mbox{if $\exn u(X)<\infty$,}
\\
0 & \mbox{otherwise.}
\end{array}
\right.
\]

By Jensen's inequality, one always has
$u (S_n/n)\le U_n/n,$ making the case $z>1$ straightforward, with
$
\lim_{n\to \infty } n^{-1}  \ln\pr ( W_{n}
   \ge z  )=\ln \pr (X=0)
$
(see Remark~\ref{rem_1}). Case $z\in (z^{**}, 1)$ admits treatment  parallel to the one we presented in Section~\ref{Sect_2}. 

It is not hard to see that the  exact asymptotics derived in Case~1 from Section~\ref{Sect_3} is still valid in this more general setting, with the obvious changes stemming from redefining the function~$f$ as
\[
f(t):= \frac{a+ t(b-a)}{u(a) + t (u(b) - u(a))}, \quad t\in [0,1].
\]
Case~2 from Section~\ref{Sect_3} could  also be extended to the general setting, all the details left to the interested reader.

\begin{figure}[ht]
\begin{center}
	\includegraphics[scale=.4]{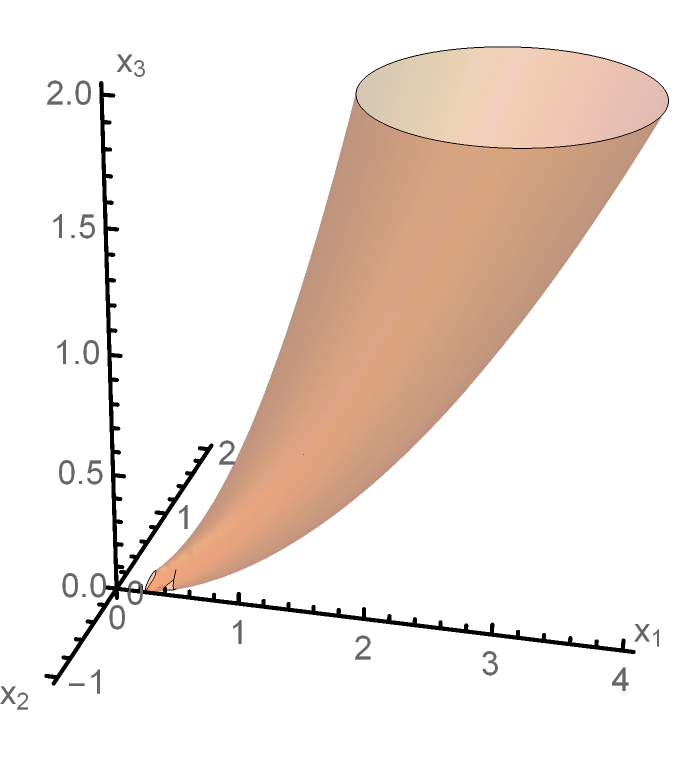}
\end{center}
	\caption{The set $B^*$ in the case when   $p=2$ and $B$ is a disk. }
	\label{Fig_3}
\end{figure}

Our approach can be further extended  to treat the multivariate case. For example, assume that $\{\bX_n=(X_{n,1}, X_{n,2})\}_{n\ge 1}$ is a sequence of   i.i.d.\ bivariate random vectors and let $\boldsymbol{S}_n:=(S_{n,1}, S_{n,2}),$ where $S_{n,i}:= \sum_{j=1}^n  X_{j,i},$ $i=1,2,$ $n\ge 1.$ For $p>1$, set
\[
T_n:=\sum_{j=1}^n (|X_{j,1}|^p + |X_{j,2}|^p),\quad n=1,2,\ldots ,
\]
and, similarly to~\eqref{Wnp}, consider the self-normalized sums
\[
\boldsymbol{W}_n:= \frac{\boldsymbol{S}_n}{n^{1-1/p}T_n^{1/p}},
\]
assuming for simplicity that $\pr (\bX=\boldsymbol{0})=0.$ Then, for a Borel set $B\subset  \R^2,$ one has
\[
\pr (\boldsymbol{W}_n \in B)
= \pr \bigg( \frac{n^{-1}\boldsymbol{S}_n}{(n^{-1} T_n )^{1/p}}\in B \bigg)
=\pr \bigg(\frac{\bZ_n}n \in B^*\bigg),
\]
where   $\bZ_n:= (S_{n,1},S_{n,2}, T_n)$ is now a three-dimensional random walk with i.i.d.\ jumps $\bz_j:= (X_{j,1}, X_{j,2},|X_{j,1}|^p + |X_{j,2}|^p)$ and $B^*:= \{\bx\in\R^3: x_3^{-1/p}(x_1,x_2)\in B, x_3>0\}.$ To illustrate this reduction, Fig.~\ref{Fig_3} shows the set $B^*$ for the disk  $B:=\{\bx=(x_1, x_2)\in \R^2 : \|\bx - (2,0.5)\|\le 0.8\}$ and $p=2$.

Evaluating the above probabilities for the self-normalized sums $\boldsymbol{W}_n$ will be an LD problem provided that, say, $\exn \bz_1\not\in B^*.$ Finding their asymptotics could be done applying the LD theory results for the three-dimensional random walk $\{\bZ_n\}.$


\begin{thebibliography}{99.}%
\bibitem{Bo13}
\kn{Borovkov, A.A.}
{Probability Theory}
{London: Springer}{2013}

\bibitem{Bo20}
\kn{Borovkov, A.A.}
{Asymptotic Analysis of Random Walks: Light-Tailed Distributions}
{Cambridge: Cambridge Univ.\ Press}{2020}


\bibitem{BoBo08}
\kn{Borovkov, A.A., Borovkov, K.A.}
{Asymptotic Analysis of Random Walks: Heavy-Tailed Distributions}
{Cambridge: Cambridge Univ.\ Press}{2008}

\bibitem{BoMo78}
\zhe{Borovkov, A.A., Mogulskii, A.A.}
{Probabilities of large deviations in topological spaces. I}
{Siberian Math.~J.}{1978}{19}{697--709}

\bibitem{BoMo92}
\zhe{Borovkov, A.A., Mogulskii, A.A.}
{Large deviations and testing of statistical hypotheses. I. Large deviations of sums of random vectors}
{Siberian Adv. Math.}{1992}{2}{52--120}

\bibitem{BoMo98}
\zhe{Borovkov, A.A., Mogulskii, A.A.}
{Integro-local limit theorems including large deviations for sums of random vectors.~I}
{Theory Probab.\ Appl.}{1998}{43}{3--17}

\bibitem{BoRo65}
\zhe{Borovkov, A.A., Rogozin, B.A.}
{On the multidimensional central limit theorem}
{Theory Probab.\ Appl.}{1965}{10}{55--62}



\bibitem{deLaSh08}
\kn{de la Pe\~na, V., Lai, T., Shao, Q.}
{Self-Normalized Processes: Limit Theory and Statistical Applications}
{Berlin: Springer}{2008}


\bibitem{DeZe98}
\kn{Dembo, A., Zeitouni, O.}
{Large Deviations Techniques and Applications. 2nd edn}
{New York: Springer}{1998}

\bibitem{Di91}
\zhe{Dinwoodie,  I.H.}
{A note on the upper bound for i.i.d. large deviations}
{Ann.\ Probab.}{1991}{19}{1732--1736}

\bibitem{GrKu89}
\zhe{Griffin, P., Kuelbs, J.}
{Self-normalized laws of the iterated logarithm}
{Ann.\ Probab.}{1989}{17}{1571--1601}

\bibitem{Mo14}
\zhe{Mogulskii, A.A.}
{On the upper bound in the large deviation principle for sums of random vectors}{Siberian Adv.\ Math.}{2014}{24}{140--152}


\bibitem{Pu01}
\kn{Puhalskii, A.}
{Large Deviations and Idempotent Probability}
{Boca Raton, FL: CRC Press}{2001}

\bibitem{Ro70}
\kn{Rockafellar, R.T.}
{Convex Analysis}
{Princeton: Princeton Univ.\ Press}{1970}


\bibitem{Sh97}
\zhe{Shao, Q.-M.}
{Self-normalized large deviations}
{Ann.\ Probab.}{1997}{25}{285--328}




\bibitem{ShWa13}
\zhe{Shao, Q.-M., Wang, Q.}
{Self-normalized limit theorems: A survey}
{Probab.\ Surveys}{2013}{10}{69--93}



\end{thebibliography}
\end{document}